# Modeling the creative process of the mind by prime numbers and a simple proof of the Riemann Hypothesis

Shi Huang, Ph.D.

The Burnham Institute for Medical Research
10901 North Torrey Pines Road
La Jolla, CA 92037, USA

Institute of Biomedical Sciences
Center for Evolutionary Biology
Fudan University
Shanghai 200433, China.

shuangtheman@yahoo.com

**Abstract**

Numbers (positive integers) are the most fundamental creatures of the human mind and the foundation to the scientific understanding of nature. Some mathematicians have suspected a link between prime numbers and secrets of creation. Understanding creativity may help resolve the deepest mysteries of primes. The algorithm that programs the mind and makes the mind creative must be sufficient for the mind to create primes. I found that primes are directly linked to the creation algorithm of the mind. The essence of primes is the duality of uniqueness and uniformity together with the creation algorithm of the mind. The creative process of the mind is lawfully determined but the outcome is unpredictable. The mathematical equivalent or model of this process is the creation of primes. Primes have the inherent property of unpredictability but can be generated by the creation algorithm of the mind, termed the Prime Law, via a fully deterministic lawful process. This new understanding of the essence of primes can deduce some of the best-known properties of primes, including the Riemann Hypothesis (RH). Understanding human creativity is obviously the most fundamental of all scientific enquiries. That this understanding can directly lead to a solution to the RH, widely considered the most important unsolved problem in mathematics, shows a deep connection between creativity and mathematics.

**Introduction:**

*Prime numbers*

A prime number is commonly defined as a positive integer that has only two divisors, 1 and itself. Both the number 1 and 2 can be either included or excluded as primes by manipulating the definition of primes. Accordingly, the primality of the number 1 and 2 are decided by human agreements rather than objective logic or reason. The number 1 is not considered a prime today but was in the past [1,2,3,4,5]. While 2 is considered a prime today, at one time it was not [6]. The odd primes have many properties not shared by 2, the only even prime. It is also easy to have a definition based on calculation that would include all primes except 2. Thus, a prime can be defined as a positive integer that cannot be expressed by the even number of sums of any single number except 1 and itself. For example, 1 is 0 (an even number) sum of 1 and itself; 3 is 2 sums of 1 and 0 sum of itself; but 2 is not a prime since it is 1 (an odd number) sum of 1.

To define numbers by calculation that is itself defined by numbers is a tautology, which merely describes ways of identifying some primes but reveals little about what a prime really is or the essence of primes. (A tautological or circular definition necessarily means a lack of true understanding.) This leads to the dilemma that 2 is a prime in one definition based on division but not a prime in another equally plausible definition based on addition. It is arbitrary human convenience to favor one tautological definition over another. We can only resolve such dilemma with objective reasoning when we achieve a deeper understanding of primes that is based on knowledge more fundamental than calculation and numbers. Primes are the foundations of mathematics and should have a form of existence or definition that is independent of mathematics.

To avoid circularity, a creature must be defined by things that are more basic than it rather than more advanced. We must use quantum particles rather than molecules to define atoms, even though we discovered molecules before we knew about quanta. Just because calculation was discovered before prime numbers in human history does not mean that primes must be defined by calculation. A concept can only be defined by concepts lower or more basic in logical hierarchy. The hierarchy must be, from bottom/basic to top/advanced, the following: numbers (primes and non-primes), addition/subtraction, multiplication/division, etc. What is even more basic than numbers must be used to define primes and non-primes. If primes are atoms that build other numbers, then the primes must be built by its own building blocks, which would be equivalent to quantum particles.

Numbers are creatures of the mind. The essence of a creature is its building blocks together with a rule of manipulating the building blocks. The essence is what is ultimately responsible for the properties of a creature. The essence of matter, the quanta building block together with a law of manipulating the quanta, is what is ultimately responsible for the properties of the physical universe. A creature must be defined by its essence. The present definition of primes contains no concept of building blocks and does not reveal the essence of primes. Without knowing the essence, it is expected that certain properties of primes simply cannot be understood.

*The creative process of the mind*

To ancient Chinese, the essence of all phenomena was the dynamic balance of bipolar opposites, yin and yang. This view was first described in the 3000 year old Chinese book *I Ching* or Yi Jing (The Laws of Changes), the most influential text in the history of Chinese civilization [7]. All opposites, day and night, happy and sad, beauty and ugly, are merely different aspects of the same phenomenon. The conflict between two interdependent opposites can never result in the total dominance of one side, but will always be a manifestation of the interplay between the two sides. Neither side can exist independent of its opposite side. In *I Ching*, the most fundamental pair of yang and yin is Heaven or Qian 乾 and Earth or Kun 坤. Qian is creative and initiative and Kun is receptive and adaptive to Qian. Nature is an endless cycling process from Qian to Kun to new Qian. The ancient Qian/yang and Kun/yin concept reveals an interdependent relationship between to create and to adapt.

In the modern interpretation of nature, what may be the most essential pair of polar opposites? The fundamental quantum particles exist in a supernatural fashion in two opposite ways, to be unique as a particle and to be a uniform part of a population wave. One of the founders of particle physics Niels Bohr considered the particle picture and the wave/pattern picture as two complementary descriptions of the same reality. Neils Bohr recognized the parallel between his concept of complementarity and the ancient Chinese thought. In choosing the *Tai-chi* symbol for his coat-of-arms, Niels Bohr acknowledged the profound harmony between ancient wisdom and modern science [8].

The creativity of the human mind is the most remarkable feature of humans that sets humans apart from all other biological species. Comparing today's civilization with those of a few thousands years ago, it is clear that humans have been constantly creating things, both physical and metaphysical. While the history of human civilization has seen a countless number of human creations, with the recent creations generally more complex than earlier ones, it seems that the basic capacities of the human brain have remained relatively unchanged at least within the last 5000 years. It seems to be a real phenomenon that more complex things get created over time while the basic capacities of the brain have stayed largely the same. The brain appears to have the ability to know or absorb whatever that have been created before and to come up with novel inventions.

There is likely a general pattern or rule that can describe the creative phenomenon of the human mind. While creations can be countless and of various kinds, the general rule or algorithm employed by the mind for each creation may be the same. Creations by successive generations overtime may be viewed as the iterative application of the same creation algorithm of the mind. To figure out what that algorithm may be is important to understand the functioning of brain cells that seem to be able to absorb past creations and to come up with something new. It may also help to reveal why the mind always becomes bored with new things after a while and why it has an insatiable appetite for novelty.

The generation of creative ideas is generally viewed as an evolutionary process. Some think it is Darwinian [9,10], while others not [11]. The question addressed by those studies is how a creative idea evolves from a population of competing ideas within a mind. The question that concerns me here is different and at a more fundamental level. It is about what motivates a mind to create in the first place. Why creation occurs and never stops in history? Why new creations get invented all the time rather than nothing gets created? Why some minds see a need to create while some other minds are content with status quo? Why faced with the same reality, some see a bottle half empty and some see half full? What is the driving force behind the mentality of a creative individual and what is the force behind those who are more content with status quo?

The creative process that I would like to analyze here is the paradigm evolution process of Kuhn [12]. A paradigm is often initiated by a creative individual and gets established subsequently by countless individuals who make incremental advances within the paradigm. For an established paradigm, some begin to see some major problems while most do not. Then a creative individual comes up with a revolutionary solution to a major problem of the established paradigm and a new paradigm is initiated. The creative process from paradigm A to B to C etc continues seemingly without end. Is the mentality responsible for paradigm change different from the mentality responsible for improving or maintaining a paradigm? Why faced with the same paradigm, some sees a need for paradigm change while others are more likely to defend it? Here, I attempt to develop a general hypothesis of the creative process from paradigm to paradigm that would apply to all human creative activities. I suggest that the human mind follows a creation algorithm.  If a computer can be programmed with this algorithm, it would become creative.

The creation algorithm that programs the mind and makes the mind creative is obviously the foundation of all human creations. Since numbers (positive integers) may represent the most fundamental creatures of the mind, the creation algorithm of the mind should be able to create numbers.  I have here found that prime numbers can model the creative process of the mind.  Others have independently noticed the connection between prime numbers and creation. The mathematician Louis Kauffamn and Hector Sabelli observed: "The generation of primes epitomizes the causal creation of novelty." [13].  The mathematician Don Zagier noted: "Upon looking at these numbers, one has the feeling of being in the presence of one of the inexplicable secrets of creation." [14].

If creation by the mind is lawful or deterministic, i.e., determined by an algorithm, it could be predictable.  However, if creation is predictable, it would no longer by definition be novel or unique.  For creation to be meaningful to humans, it must not be predictable.  Human creativity is also logical and reasoned and does not seem to be arbitrary.  Thus, whatever the algorithm that has programmed the human mind and made the mind creative must make the creative process lawful and yet the outcome unpredictable.  Similarly, prime numbers have been found to be both lawful and seemingly unpredictable. But a deterministic law of primes or of human creativity remains to be discovered that would make the process lawful but nonetheless allow the outcomes unpredictable.

*The Riemann Hypothesis*

The German mathematician Bernhard Riemann formulated the Riemann Hypothesis (RH) in 1859 [15].  The hypothesis is widely regarded as the most important unsolved problem in all of mathematics. In an interview, the mathematician David Hilbert explained that he believed the RH to be the most important problem 'not only in mathematics but absolutely the most important.' [16].  The Fields medalist Enrico Bombieri said in an interview: "The Riemann Hypothesis is not just a problem.  It is *the* problem.  It is the most important problem in pure mathematics.  It's an indication of something extremely deep and fundamental that we cannot grasp." [17].  The RH is believed by most mathematicians to be true.  A large number of deep and important other results have been proven under the condition that it holds.

The RH essentially says that the primes are as regularly distributed as possible given their seemingly unpredictable occurrence on the number line.  According to the Prime Number Theorem of Gauss, the number of primes less than N is approximately the logarithmic integral

Li(N) or less precisely N/ln(N). If the RH is true, the error between Li(N) and the true number of primes is at most of the order of the square root of N [16,18,19,20,21]. This error margin is the smallest possible and cannot be improved by much [22]. This is the error margin expected by the theory of probability for some unpredictable events such as a coin toss. If one can prove that primes must have a regularity pattern which must have an error margin which must be the smallest possible, one proves the RH.

Primes are both lawful and seemingly unpredictable, as suggested by the RH. This is highly similar to creations of the mind. If the RH is a truly fundamental problem of primes and if primes are truly related to secrets of creation, as some have suspected, then the RH should simply remain unsolvable until we have some understanding of creativity. Here, I first describe the creation algorithm that makes the creative process lawful but the outcome unpredictable. I then show that this algorithm can create primes in a fully deterministic and lawful process but still allows the primes to have the intrinsic property of unpredictability. By showing that primes are lawfully generated but unpredictable, I prove that primes must have a regularity pattern which must have an error margin which must be the smallest possible. This in effect proves the RH. I further prove using simple logic that it is impossible for conventional mathematics to prove the RH. Discovering the deterministic law of primes or of human creativity is the only possible way to prove the RH.

**Results and discussion:**

**The algorithm that makes the mind creative**

The most fundamental capacities of a human mind may be to know and to imagine, which are essential to creativity. To know is to recognize the unique from a background of contrasting uniformity and vice versa. To imagine is to think of novel things that do not exist previously. By observing how the human mind creates, I have found a simple algorithm that programs the mind and makes the mind creative. This algorithm consists of a pair of opposite but complimentary yin and yang principles and a mind that coordinates the interplay of the two principles. A creation or creature is defined as the unique that does not exist previously, is distinguished from all other imagined things, and can exist subsequent to its creation by being able to initiate a population of followers that share a uniform pattern resembling the unique. A creation has the bipolar duality of uniqueness and uniformity. A follower of a creation is defined as the new thing that does not exist previously but shares some uniform property with a prior creation. A creation is a big jump in paradigm while a follower of a creation represents a small step advance within a paradigm.

The imagination of a mind is either following the existing patterns of past creations or is based on a novel pattern. How does a novel pattern come to the mind remains a mystery and is of no concern here. A new but meaningless thing or pattern is not a creation because it cannot be uniquely distinguished from other imagined entities, cannot be logically linked with existing patterns, and cannot initiate a following. A great piece of music or book or art initiates a following by existing in the minds of other people. A book that was soon forgotten forever is not a creation but is merely a follower of an existing pattern. Existing pattern consists of both past creations and of a default order-less state. The order-less state is the background and driving force for order and pattern. A new thing that is not following any existing ordered pattern but is not uniquely distinguished from the order-less state is still viewed as a follower because it is following the existing order-less pattern. Things that constitute the order-less state include all that cannot be logically linked to any ordered pattern and cannot be uniquely distinguished from others or are equally unique as others.

The yang principle is uniformity selection that allows the mind to recognize the unique or the creation. Uniformity abolishes individuality and selects for the unique. Uniformity selection drives the creation of the unique. The yin principle is uniqueness selection that allows the unique to initiate a population of followers sharing a uniformity pattern resembling the unique. The mind uses this principle to allow the unique to exist or survive subsequent to its creation. Uniqueness selection results in the formation of an ordered uniformity consisting of individuals that are fittest or most adapted to the unique. The process from the unique to a specific uniformity of a population of followers is essential for the unique to exist subsequent to its creation, which further serves to drive the creation of the next unique. The creative process of the mind is the iterative use of the same creation algorithm and an endless cycling process from uniformity to unique to new-uniformity. When the mind sees the unique, the mind strives to fit and follow. When the mind sees uniformity, the mind strives to be unique. All human minds are a unity of different degrees of the yin and yang principles.

To create, the mind needs to know what is known previously, which is termed the existing-uniformity. Selection by existing-uniformity allows the mind to know whether something is new with a meaningfully ordered pattern. In addition, all creations begin from the imagination of the mind. Within the imagined world, there exists a unique entity that is distinguished from the imagined-uniformity shared by other imagined entities. To create by uniformity selection is to bring into existence an imagined entity that is distinct from both the existing-uniformity and the imagined-uniformity.

The formation of the order-less existing uniformity is by the default of reproduction and the inherent nature of the mind. The mind treats anything that cannot be rationalized with past creations as part of an order-less uniformity. Random brushes on a canvas would belong to the order-less uniformity. The formation of the ordered existing-uniformity requires the yin principle of uniqueness selection. This selection process selects individuals to follow the unique creation of the past. The followers of a unique creation are essential to the popularization of the unique and the long-term existence of the unique in the form of existing-uniformity. The followers also contribute new variations or incremental advances around the main theme/paradigm of the unique creation, which would form a new level of existing-uniformity essential for triggering the next unique creation. However, the incremental progress made by the followers cannot directly in itself lead to the next unique creation. Creation of the unique represents a discontinuous change in paradigm as defined by Thomas Kuhn and is fundamentally different from the formation of followers.

Existing uniformity thus consists of order-less or yin and ordered or yang. Based on the existing uniformity, the mind is able to know whether something is imagined or not yet existing. Among the things imagined, a uniform property may be shared by all except the unique. The unique is the one that has the closest relationship to the existing uniformity but does not belong to any of the existing paradigms. As noted by the author de Bono: "Every valuable creative idea must always be logical in hindsight. If an idea were not logical in hindsight then we would have no way of seeing the value of that idea and it would simply be a crazy idea."[23]. The creation of the unique cannot come as a logical extension of an existing pattern but is nevertheless logically related to existing patterns after the fact of creation.

**Creating primes by the creation algorithm of the mind**

Like creations of the mind, the odd primes including the number 1 also have the dual property of uniqueness and uniformity. A thing is unique if it is not an inherent part of something else and is different from uniformity. A number is an inherent part of a smaller number either because it is needed for the smaller number to have meaning or because it can be expressed as a pattern of a single smaller number >1. The number 2 lacks uniqueness because it is an inherent part of creating the number 1, as evidenced by the existence of civilizations that had invented only 1 and 2 and by the absence of civilizations that invented only 1 but not 2. We need 2 to invent 1 or for 1 to have any meaning. We need both 1 and 2 in order to invent the concept of number.  However, we do not need 3 to invent 1 and 2.  The author John Barrow observes:

"There are many primitive human groups who cannot count beyond two and have developed number-sense at all.  Some Australian aboriginal tribes only possess words for the quantities 'one' and 'two'.  All greater quantities are expressed by a word with the sense of 'many'.  Some South African tribes exhibit a similar language structure for number words.  In fact, there is even a vestigial remnant of it in modern Indo-European languages where we find that the original root for 'three' has a meaning of something like 'over' 'beyond', or 'afar', indicating that the early sense was merely that of something more than the particular words for one and two.  In Latin, this affinity exists between *trans* meaning 'beyond' and *tres* meaning 'three'; the same relationship is manifest in French where we find *tres* for 'very' and *trios* for the number 'three'.  Presumably, only later did words emerge to distinguish the various different varieties of 'many'.  Some languages, for example Arabic, also retain a threefold treatment of quantity with differentiation between singular, dual (for two only), and plural (for more than two)." [24].

All numbers are inherent in the number 1 as patterns of 1s but the property of uniqueness of the odd primes is not inherent in the pattern of 1s.  Uniqueness is in contrast to uniformity and cannot exist independent of uniformity. While a prime can be expressed as a pattern of 1s, its uniqueness cannot.  Every number (positive integer) can be uniquely defined by a pattern of 1s but this makes every number equally unique.  Thus none is unique.  The uniqueness of a number is based on the existence of numbers greater than 1 and the existence of non-unique numbers.  Primes and non-primes are like odd and even or yang and yin and cannot exist independent of each other.  The number 1 is unique since oneness is synonymous with uniqueness. If 1 is unique, then 2 must be non-unique because it is an inherent part of creating the number 1.

The opposite of uniqueness is uniformity or not being able to be singled out.  A prime also exists in a pattern, e.g., 18 is a pattern of the prime number 3.  In such a pattern, the number 3 could not stand out as a unique individual. The uniformity property of a prime makes it possible for other subsequent primes to be uncovered as the unique.  The number 23 is a prime because it is not a pattern of any other numbers greater than 1.  The number 2 is essential for the number 1 to be unique and for other odd primes to be unique.  For uniqueness to exist, the uniformity background must co-exist.  Two is the first number of non-uniqueness and therefore has some uniqueness property and the related uniformity property.  It is the most unique (the first number of non-uniqueness) and the most uniform (present in more patterns of 1s than any other non-unique numbers) among non-unique numbers.

If the building block of non-primes is the prime, it is only fair and logical to go down the hierarchy to ask what is the building block of prime. That building block logically cannot be a number since prime number is the lowest level a number (positive integer) can be. If 1 is a prime, its building block must be 1 itself. The number 1 is also the building block of all other primes. A prime is a positive integer that can be built in only one way from its building block 1 by way of even number of sums of 1 but not of any other numbers greater than 1. How does a creative mind perceive the number 1? Of course, 1 represents uniqueness or oneness or a single smallest point of the whole. One is also uniformity or the single wholeness and is present everywhere or in every number or in every part of the whole. So, 1 embodies the ultimate duality of uniqueness and uniformity. To a creative mind, 1 and the duality are synonymous. The number equivalent of the duality concept is 1. Since 1 is the sole building block of primes, we can also say that the duality is the building block of primes. The duality as building block of primes expresses the meaning of 1 as building block in a more fundamental way that is directly linked to the creation algorithm of the mind. It can also be easily related to the quantum building block of matter, which also has the duality of uniqueness/particle and uniformity/wave. The following shows that the creation algorithm can use the duality as building block to create primes. The mathematical model of the creation algorithm is the orderly creation of primes.

Postulate 1. The imagined domain  All things created by the mind comes from imagination and the imagined world of the mind is termed the imagined domain. The content of this domain consists of an infinite number of the basic building block of numbers, 1. There are infinite number of patterns of 1, each differ by its count of 1s. Each pattern, except that of a single 1, has the uniform property of having a count of 1s that is between two other patterns. The pattern of 2 is between the pattern of 1 and the pattern of 3. Since the contents of the imagined domain has no numbers smaller than 1, the pattern of a single 1 is not in between two other patterns and is therefore unique.

Postulate 2.  The reality domain  The reality domain is where the materialized creations of the mind exist. A prime is generated in the reality domain because of its uniqueness at the time of its creation. It subsequently exists in the reality domain because of its ability to initiate a pattern/uniformity. A prime is defined as a lawful creature of the mind that has the duality of uniqueness and uniformity. A non-prime is defined as a follower of a prime. The mind creates primes by following the two principles of the creation algorithm as postulated above: 1) to generate uniqueness by uniformity selection and 2) to maintain subsequent existence of the unique by uniqueness selection to form uniformity. Uniqueness selection is the process of species formation or forming follower numbers that share properties with the unique. For example, the follower numbers of 3 are 6, 9, 12, . . ., which share the uniform property of 3-ness and form the species of 3. A pattern of 1s or number moves from the imagined domain into the reality domain because it is either uniquely recognized by the mind or is necessary to maintain existence of the unique in the reality domain.

<u>Creating primes.</u>  Prior to the creation of any numbers in the reality domain, the unique number in the imagined domain is 1. So the first goal is to generate 1 as the unique or prime in the reality domain. Since a prime must form a species or pattern in order to exist following its creation, the species of 1 is formed with 1 followed by the next closest number 2. In addition, to express uniqueness requires the simultaneous presence of uniformity. So, the species of 2 is formed to represent uniformity with 2 followed by the next closest number that shares the property of 2-ness, 4. Two is selected to represent uniformity because it is the only other pattern besides 1 that is available in the reality domain at this point when the species of 1 has not yet progressed beyond 2. Table 1A shows the contents of the reality domain at its time of creation. The prime/uniqueness/1/odd/yang and non-prime/uniformity/2/even/yin are generated simultaneously and cannot exist independent of each other.

After the beginning stage of generating the reality domain, the mind is aware of both the imagined domain and the reality domain. By comparing the two domains, the mind is looking for the next prime or unique pattern among patterns in the imagined domain that have no match in the reality domain. This pattern is now 3 and it is unique because it is the smallest while all other patterns share the uniform property of having counts of 1s that are between two patterns. To express 3 as a prime, the species of 3 (3, 6, 9) is formed in the reality domain. To apply the new concept of 3-ness, all species are extended to the $3^{rd}$ position. The reality domain has now advanced from the beginning stage of 1 and 2 to the next stage of 3-ness. At this stage, a number larger than 3 such as 4 expresses only the concepts already established such as 2-ness (2 units of 2-ness). After the stage of 3-ness, the mind is again ready to look for the next unique pattern remaining in the imagined domain, which is now 5. From the concepts of 3-ness and 5-ness, the 4-ness of 4 is now recognized as the intermediate between 3 and 5. By applying the concept of 5-ness and 4-ness, all number species are extended to the $5^{th}$ position. The species of 5 (5, 10, 15, 20, 25) is formed to express 5 as a prime. The mind is then ready to look for the next unique pattern that remains in the imagined domain (Table 1B). In this way of iteratively applying the same creation algorithm, an infinite number of primes can be generated. Because this creation algorithm of the mind can create primes, it is hereafter called the Prime Law. Since primes have the same property and meaning as creations of the mind, the word 'prime' and the word 'creation' are interchangeable or synonymous. Therefore, the 'Prime' Law also literally means the 'Creation' Law.

**Discussions on creativity**

*The yang and yin principles of the Prime Law*

The uniformity selection principle suggests that the mind is capable of converting all that exist into a background upon which to base new imaginations. The uniqueness selection principle suggests that the mind is also programmed to adapt to existing paradigms. A human mind feels the need to fit in with the conventions of society but also feels the need to be unique or different from all other people. Humans display polar opposite sides of creation-related character traits that are selected to coexist by the yang and yin principles. The uniformity selection principle values individualism, ambition, adventurism, self-centeredness, and distaste for routine labor, while the opposites are valued by the uniqueness selection principle. Both sides are essential for creation to go on, and all humans display unity of different degrees of both. The uniqueness selection principle also accounts for the inherent drive of humans to publicize their creative work once they have created something. If they do not work hard to present their creative work to the public and to have their work accepted and followed by others, their work would not count as a complete creation and would have no impact on the creative process of humanity. The Prime Law suggests that this drive to have others to accept and follow one's own creative work is hardwired into the brain and is essential to human creativity.

The Prime Law is consistent with the ancient Qian and Kun concept and further develops that concept. The ancient made it clear that Kun is to adapt to Qian so that the initiation of Qian can be brought to completion. But they did not say that Qian or creativity is dependent on Kun. They also did not express the concept of uniqueness and uniformity and the principle of selection. The Qian and Kun concept does not invoke the mind and has not previously been used to account for the creative algorithm of the mind. The uniformity selection principle suggested here is a novel concept. The uniformity formed by the adaptive process of Kun forms the background and inspiration for the creative imaginations of Qian. Kun needs Qian so it has something to adapt to. Qian needs Kun for Qian's creation to take form and to exist. Qian also needs the products of the Kun process as the necessary background for new creations.

The novel concept of uniformity selection seems obvious and explains the mind's insatiable appetite for novelty. The idea was inspired by death selection in nature or the inherent drive of life to stay away from death. Death is not disappearance of matter but is merely a return of matter from a unique high complexity ordered state (we call life) to a uniformity state of less-ordered matter (we call death). All life becomes the same in death in terms of matter. Death is order-less uniformity state of life-building molecules. Each human mind has an inherent need to know what is the self or me or what is special or unique about the self. For the mind to stay away from uniformity/death, the mind needs to be unique. The only way to be unique is to be creative. But in order to be unique and creative, the mind needs to know or learn first what is uniformity.

Uniformity selection describes the big jump creations and uniqueness selection describes the small step creations. Small step creations are creations within a paradigm. Big jump creations are changes in paradigm. It is widely noted that small step creations cannot add up to big jump creations. As the author Edward de Bono notes: "It must be said that a succession of small jumps do not add up to a big jump. A big jump is usually a paradigm change or a new concept. Because this can involve a total reorganization of previous concepts, this is not likely to come from an accumulation of small jumps." [23].

The Darwinian theory of evolution is a creation law by a mind-less process. Mind is not needed in such a law for creation to occur. It is therefore hard to imagine that a mindless creation law can accurately and completely describe the creation process of the mind. We would expect that the creation law employed by the mind would be different from that employed by a mindless process. Nonetheless, the creative process of the mind has been viewed as a Darwinian process of blind-variation and selective retention [9,10]. Thought is considered as a series of tiny selections. Thus one arrives at an idea by, given one's current thought, selecting among the possible ways of varying it to generate a next thought, and repeating this process until the idea appears. The generation of ideas is blind without foreknowledge of what it will eventually lead to or result in.

While the notion of blind generation of ideas is debatable [11], this Darwinian view is actually a good description of the process of incremental advances within a paradigm. For example, within the Darwinian paradigm, various sub-branches of Darwinism such as the Darwinian view of the creation process emerged as slightly new ideas. The creator of a sub-branch of Darwinism generates new variations of Darwinism to solve a new problem and selects the one that appears to solve the problem but also fits best with the Darwinian paradigm. That creator is acting via the uniqueness selection principle. He grants the Darwinian paradigm to be true and generates variations of it and selects the one that fits with Darwinism while also solves a new problem. He is using the existing paradigms to solve new problems. His driving force is to follow and improve the existing paradigms. Under such a driving force, he would have no chance of creating an idea that would expose fundamental flaws with the existing paradigm. His efforts will only serve to strength the existing paradigm. He would not have developed an evolution theory of the creative process if he has not known Darwinism in the first place. The generation of his ideas are not random and blind but are tinted by the existing paradigm. The selection of his best idea is also strongly influenced by how fit it is to the existing paradigm.

In the process of applying Darwinism to new fields, it may be discovered that Darwinism may not apply well in certain fields, which may inspire the creation of a better and substantially new theory. The application of Darwinism by its followers creates a new level of uniformity that serves to drive the next unique creation. But so long the creator is driven to generate and select his ideas based on Darwinism, he is unlikely to create anything that differs from Darwinism.

In contrast, a person who is motivated by the principle of uniformity selection is less likely to follow any existing paradigms and is more motivated to invent a new paradigm. Such a person would have little interest in working in a field where the chance for paradigm shift is minimal. His ideas may come from fields that are not directly related to the field of his main interest. But when he connects those ideas from other field with a new field, he sees a unique combination of idea and problem. His selection criteria for the best idea are 1) if it is substantially new as a solution to a problem, and 2) if it is unique among all his ideas in that it can be best linked logically with all relevant existing paradigms despite not being an obvious variation of any existing paradigms. Here, the critical idea does not originate from a variation of a preexisting theme. The idea may come from a different field or from novel insights in studying a different field. The idea may not be a logical extension or variation of any existing paradigms but may be logically linked with them in hindsight. No variations from realism could have led Picasso to invent cubism. No variations from the Darwinian view of the creative process could lead to the uniformity selection principle as proposed here.

*Purpose of the mind*

Since the arrival of humans, countless creations have been achieved by the human mind. Will the mind run out of things to create? The Prime law suggests that the mind will never run out of things to create. Creations are primes. Primes are infinite and so are creations. Why is there something created rather than nothing? This is perhaps because the mind has an inherent drive and capacity to imagine new things based on existing uniformity. The mind may have an inherent purpose of endlessly creating uniqueness. The purpose of the mind is to create. This purpose can be achieved by the iterative application of the creation algorithm or the Prime Law as proposed here. Indeed, for the human mind to have a purpose, that purpose must be achievable and yet endless.  If it is not endless, one could ask what would be the next purpose when a prior purpose is achieved.  To create uniqueness or primes is achievable and yet endless, as demonstrated here using prime numbers as models of creations.  Each prime is a unique creation but the number of primes or creations can be proven to be infinite.

*Creation is not merely something new*

To create (uniqueness) is more than simply making something new.  To be new is necessary but not sufficient for a creation.  A creation is not only new but is also unique among all potential new things that can be imagined by the mind and its uniqueness lies in the fact that it has the closest relationship to, or is the most smoothly adapted to, the whole pattern of what exist previously.  This concept can be illustrated by the prime generation scheme as shown in Table 1B, which represents a time of existence that has no concept beyond 5-ness.  Many numbers are missing in Table 1B and can qualify as new, such as 7, 11, 13, 14, 16, 17, 18, 19, 21, 22, 23, 24, and all numbers larger than 25.  But only 7 is unique because it is the smallest missing number or because all numbers smaller than the unique number belong to what exist previously. A prime number or a creation is one that is best adapted to the whole pattern of all that exist, whereas a non-prime number or something that is merely new is best adapted to a particular sub-pattern of the whole.

Why creation has the property of uniqueness?  Because if it were not, the mind would not be able to recognize it from an infinity of choices.  The mind does not create or select ideas by throwing a dice.  If that were the case, the creations would not have the property of uniqueness, and the mind would be unconscious of the properties of such creatures.  Such creatures would lack any coherent logical relationship among them and would not be able to form the uniformity pattern to drive the next creation.  To create by throwing a dice is meaningless to the mind.

*Discovering laws is creation*

Mathematics has been viewed by some as the creation of the mind.  That the essence of mathematics, the primes, has properties of the creations of the mind is remarkable. Even if mathematics exist in a Platonic world independent of humans, the discovery of mathematics by humans involves the same creative process that human employs for creation. To create means to bring something into existence within the world of human minds.  That something could be a natural law or mathematical law that has existed independent of humans or it could be a painting that has never existed anywhere prior to its creation. A scientist only discovers laws but the process of discovery is no different from the process of creation, in the sense that something unknown is now made known to the mind.  However, the discovery of laws involves much more creativity than the chance discovery of physical phenomena or facts.  A false law is clearly a human creation.  Since a scientific law or theory or hypothesis has the potential to be false or falsified, it is always a human creation regardless false or true.

*Testing the Prime Law*

Can the Prime Law be tested or falsified?  The uniqueness selection principle is essentially Darwinian and can be easily verified by observing how students learn and follow rules or past knowledge.  Students who are the fittest or the best in learning past knowledge are following the unique of the past or are more driven by the uniqueness selection principle.  The Prime Law predicts that such students are less likely to be creative or be driven by the uniformity selection principle.  They are less likely to be bored with just learning.  The Prime Law predicts that the creative person should be someone who is fit but far from the fittest.  To test this, we could study the history of the past creators who have made outstanding intellectual contributions to human culture.  For example, was Confusions, Socrates, Picasso, or Einstein fittest at the time when he made his creations?  We could examine how many creative works in history were done by the fittest who had power and money at their time.  On the other hand, the Prime Law also predicts that a person who cannot learn or become somewhat fit due to either economic or physiological reasons are not likely to be creative.  We could examine whether people with learning difficulties have a lower chance of being creative.

*How to be creative?*

The creation algorithm developed here explains how and why minds have created countless creatures in history and will continue to create endlessly.  It also provides practical guidance for a human being to become creative. Each person is a unity of the yin and yang principles and desires to be unique as well as to fit.  To create is at a minimum to be unique and vice versa.  The yin principle of uniqueness selection directs the growth and learning period of a budding creator and the yang applies in the mature period of the creator. By not following the yin principle of uniqueness selection, the budding creator will remain a part of the default order-less existing uniformity, or may keep his/her creative work private rather than working hard to present the work to others in order for others to accept and follow the work.  By following only the yin but not the yang principle, the creator will become a member of the ordered existing uniformity or a fittest member of an existing paradigm. By following only the yang principle, the creator's work may never be recognized by society. Good works are products of the yin principle whereas great works are products of unity of the yin and yang principles.

**On indivisibility of primes**

Uniqueness means that a number is not an inherent part of a smaller number greater than 1. A prime is not a pattern of any smaller number greater than 1, which means indivisible by any smaller number greater than 1. Indivisibility is therefore a secondary property of primes as the unique and should not in and of itself confer primality. The number 2 is indivisible but is not a prime because it lacks the uniqueness essence. It is an inherent part of creating the number 1 as the unique or prime.

**On the duality of unpredictability and regularity of primes**

It is well known that primes seem to exhibit the duality of unpredictability and regularity. Such seemingly impossible unity of extreme bipolar yin-yang opposites is what makes primes so interesting and mysterious. The mathematician Don Zagier once said in a 1975 lecture: "There are two facts about the distribution of prime numbers of which I hope to convince you so overwhelmingly that they will be permanently engraved in your hearts. The first is that, despite their simple definition and role as the building blocks of the natural numbers, the prime numbers grow like weeds among the natural numbers, seeming to obey no other law than that of chance, and nobody can predict where the next one will sprout. The second fact is even more astonishing, for it states just the opposite: that the prime numbers exhibit stunning regularity, that there are laws governing their behavior, and that they obey these laws with almost military precision." [14]. However, the first fact of unpredictability remains unproven. The mathematician Robert C. Vaughan said that "It is evident that the primes are randomly distributed but, unfortunately, we don't know what 'random' means." It is the seeming randomness or unpredictability that makes the regularity of primes so striking and interesting. If primes were ever proven to be predictable, it would cease to have the duality property and would in turn cease to be interesting and mysterious.

A proof on the unpredictability of primes would also have practical benefits since it would save people from wasting time on trying to develop formulas to predict primes. There exist a variety of formulas for either producing the N th prime as a function of N or taking on only prime values. However, all such formulas require either extremely accurate knowledge of some unknown constant, or effectively require knowledge of the primes ahead of time in order to use the formula [25]. They do not really count as prediction. A true predictive formula should not make use of the knowledge of existing primes in order to predict the next future prime. It must be the opposite of the non-predictive way of generating primes, which is to use knowledge of existing primes such as indivisible by known primes.

Since the apparent duality of unpredictability and regularity is what makes the primes interesting, it is essential to prove that the duality is real. The newly discovered essence of primes can easily deduce or prove the duality of unpredictability and regularity. Primes can be easily proven to be unpredictable given the new essence of primes as described here. When something can be predicted, it must belong to a pattern. As such, it is not unique and hence, by definition, not a prime. The essence of uniqueness rules out prediction of primes as a viable possibility. There is also another easy way to prove this. To predict primes means to predict uniqueness and in turn uniformity since uniqueness needs uniformity to have meaning. Uniformity is made of existing primes. So to predict primes is to predict *existing* primes, which is a logical non-sense.

The uniformity essence of primes demands that the formation of uniformity from a newly created prime or uniqueness must be regular and predictable. So the uniformity forming property of primes gives rise to the regularity of primes. Primes exist as regularly as possible in the uniformity. The creation of primes is fully determined by the orderly formation of uniformity by existing primes. The lawful rather than lawless way of creating primes explains why primes should follow some regularity patterns, such as the Prime Number Theorem. The unpredictability of individual primes explains why such a pattern cannot be completely precise or free of error margins. Thus, the essence of primes fully explains the duality of unpredictability and regularity.

**On proving the RH**
*Seeming randomness and real randomness*

I use the phrase 'seeming randomness or deterministic randomness' to describe an outcome of a *lawful* process that is nonetheless unpredictable, like the creation of primes by the Prime Law. A population of such seemingly random outcomes should show a regularity pattern reflecting the lawfulness and regularity in the process leading to these outcomes. However, even the most precise pattern should still show some error margin reflecting the unpredictability or seeming randomness of the individual outcome. Primes have been found by many to show 'deterministic randomness' [26,27,28,29,30].

I define 'real randomness' as an unpredictable outcome of a lawless/arbitrary process like selecting a prime number from an infinity of numbers by playing a dice. Here the dice throw per se is not lawless/arbitrary/random. The lawless/arbitrary/random component in a lawless process involving the dice is connecting the dice arbitrarily with a meaningful concept or event that has no lawful connection to the dice, such as connecting prime numbers with the landing of a dice. For a process that involves both a lawful component (dice throw per se) and a lawless component (arbitrarily linking landing of dice with calling a number prime), the process is effectively lawless/arbitrary/random.

Both seeming-randomness and real-randomness are unpredictable but the error margins from a pattern are greater with real-randomness. A population of lawfully caused and predictable outcomes follows a precise pattern without any error margin. A population of lawfully caused but unpredictable outcomes follows a less precise pattern with some error margin like the square root of N. A population of lawlessly caused outcomes follows a rough pattern with huge error margins which could be so high as to render the pattern meaningless or equivalent to no pattern. If whatever number that is selected from an infinity of numbers by playing a dice is defined as primes, we would obviously detect no meaningful patterns of primes in most cases, which is equivalent to saying that we could only have patterns with huge error margins. The error margin for a pattern of outcomes that are lawfully caused but unpredictable must necessarily be the smallest among patterns that cannot predict individual outcomes. Any smaller error margin would mean some degree of predictability. If we know that certain position of the tossing hand could cause a higher chance of landing heads while another position favoring tails, we could improve on the error margin but then the coin toss would not qualify as truly unpredictable. True unpredictability is shared by all kinds of randomness. Among these, the seeming-randomness or unpredictability of outcomes of a fully lawful process has the least amount of randomness or the smallest error margin from a regularity pattern.

There is a pattern that a fair coin toss follows, which says that the number of heads is equal to half of the number of toss N with an error of the square root of N. This pattern is a law that is valid based on logical reasoning alone. A fair coin toss must not have irregular or arbitrary/random bias toward the head or tail. Each landing of head or tail is fully determined by laws, such as the gravitational law, the exact position of the tossing hand, the wind, etc. A lawful process should produce reproducible outcomes. A coin toss is reproducible if the tossing conditions can be exactly reproduced. A coin toss is only seemingly random because of unpredictability. It is unpredictable because humans cannot measure all the physical parameters that determine the fall of a coin. Also, the laws are not biased to favor either head or tail and remain unchanged timelessly. If a divine were to suddenly intervene for no reason to cause more landing of the head, the coin toss would be lawlessly caused and would display much wider error margin. If we only detect a seeming randomness in our coin toss with an error of the square root of N, we would be confident that everything is well and regular and no laws have been broken by either random accidents or deliberate intentions. But if we see a much wider variation than the square root of N, we would know that something is wrong or that some laws have been broken either accidentally or deliberately. The coin toss would be considered as unfair.

*A simple proof of the RH*

Since the RH means that the error margin from the pattern Li(N) is the smallest possible, one can prove the RH by showing that primes must have a regularity pattern which must have an error margin which must be the smallest possible. I have proven this in principle by showing that primes are lawfully caused and yet unpredictable. The creation process of prime is lawful and non-random or non-arbitrary. But the outcome of this process, i.e., calling a number a prime, is unpredictable. Because of the unpredictability, one simply cannot have a pattern of primes that is free of error margins. (A pattern without error margins would mean predictability.) However, because of the lawfulness of the process, the error margins must be the smallest possible among all kinds of unpredictable outcomes, which include those caused by either lawfully or lawlessly determined processes. Specifically, it must be smaller than the error margin of outcomes that involve a lawless process such as arbitrarily calling 6 a prime or calling the head of coin prime. The lawful creation of primes is similar to the phenomenon of coin toss. It is therefore expected that the two phenomena should have similar error margins. In both cases, the error margins are the smallest possible. Any bigger error margin would mean some degree of lawlessness in the process of creating primes or in the process of coin toss.

The mathematician Harold M. Edwards wrote: "One of the things which makes the Riemann Hypothesis so difficult is the fact that there is no plausibility argument, no hint of a reason, however, un-rigorous, why it should be true." [18]. The essence of primes is the plausibility reason that the RH must be true. Given this essence, is it possible for the RH to be false? The falsity of the RH would indicate arbitrariness or lawlessness in the process of creating primes. According to Fields medalist Enrico Bombieri, "The failure of the Riemann hypothesis would create havoc in the distribution of prime numbers."[22]. But there is nothing unlawful or disorderly in the process of creating primes by the Prime Law. Thus, it is impossible for the RH to be false. It can only be true.

All other methods of generating or finding primes, such as division by smaller numbers and the sieve of Eratosthenes, are also orderly or deterministic. But unlike the Prime Law here, these methods cannot prove the unpredictability of primes. They define primes by the process of generating primes and thus do not give primes any meaning that is independent of the process: primes are whatever that are found by the process. Under the Prime Law, however, primes have the meaning of uniqueness/uniformity that is independent of the process of creating primes by the Prime Law. Uniqueness has meanings that are independent of the process of creating uniqueness. A creature has meanings that are independent of the process of creating the creature. Unpredictability is a property of the outcome and is independent of the process leading to the outcome. Both lawful or lawless processes can lead to unpredictable outcomes. That a coin falls half of the time head is an inherent property of the coin and is independent of the process of coin toss. If the property of the outcome is all defined or given by the process, then a lawful process simply cannot give the outcome the property of unpredictability or seeming randomness.

*Simple proofs for the impossibility of proving the RH by conventional mathematics*

Some mathematicians have suggested that Gauss had modeled his Prime Number Theorem by tossing a coin [16]. To Gauss, primes look like they are generated by nature flipping a coin, heads it's prime, tails it's not. The coin could be weighted so that instead of landing heads half the time, it lands heads with probability of 1/ln(N). If the prime coin is tossed in an unbiased and unpredictable fashion, it should show an error margin of the square root of N, the same as a real coin toss. The RH says that the prime coin is fair and lawful with an error margin similar to that of a real coin toss.

Obviously, if primes were generated by a prime coin toss, they would not be reproducible. If nature were to toss it again, 5 may not show up as heads or primes while it did the first time, even though the number of heads or primes under a number N would stay similar or reproducible. However, primes generated by the Prime Law are rigid and precise and reproducible every time. Also, if primes were generated by a prime coin toss, one would expect to find at least 2 heads in a row at some point on the number line. But odd primes are always separated by at least one non-prime number. Also, the prime coin toss would not be biased toward only odd numbers. So, individual prime is clearly not generated in an arbitrary/lawless/random fashion by nature flipping a prime coin, which involves an arbitrary/lawless/random connection between primes and landing of heads, even though the regularity pattern of a population of primes can be modeled by a coin toss. The probability for a prime coin toss to generate the actual prime sequence is as small as anyone can imagine although not zero. It reaches zero asymptotically. It is effectively zero for all practical and sensible purposes.

Mathematicians commonly try to use conventional mathematics to prove literally the exact statement of the RH, all non-trivial zeros of the zeta function have real part one-half. One can easily prove in several different ways that it is not possible for such a mathematical proof of the RH. A conventional mathematical proof of the RH is self-defeating. It will tell us nothing about the process of creating primes. Thus any process is admissible or possible, including a lawless/random/arbitrary process such as naming primes by a coin toss. Indeed many have speculated that, if the RH is true, primes would look like being created by nature flipping a coin, as the mathematician Marcus du Sautoy put it: "Physicists have grown used to the idea that a quantum dice decides the fate of the universe, randomly choosing at each throw where scientists will find matter. But it is something of an embarrassment to have to admit that these fundamental numbers (primes) on which mathematics is based appear to have been laid out by Nature flipping a coin, deciding at each toss the fate of each number. Randomness and chaos are anathema to the mathematician."[16].

If the RH is proved without knowing anything about the process of creating primes, it would mean that the process could be either lawful or lawless. We simply don't know either way. One could be perfectly justified in saying that it is possible for a lawless process to create extreme regularity even though the probability for that is essentially zero or asymptotically zero. The notion that it is possible for an effectively zero probability event to actually happen is against nature and common sense logic. No sane person can believe that a coin toss can have any realistic chance of generating the correct prime sequence under a large enough number such as numbers greater than 1000. Common sense also demands that it is impossible for an arbitrary process or coin toss to always produce the desired result for an infinite number of times.

Common sense or intuition or axioms are the foundations of mathematics. It is not possible to imagine that the solution to the most fundamental problem of mathematics should nullify the value of common sense or the foundation of mathematics. If the common sense is nullified, any deductions from it are of course also no good, which would include the mathematical proof of the RH. Thus proving the RH in the mathematical way is self-defeating. If it is done, it would invalidate the whole enterprise of mathematics based on common sense and reason/logic. To prove the RH by conventional mathematics is to prove the impossible (an asymptotical zero probability event) possible, which is logically non-sensible. Thus, by the classic mathematical method of *reductio ad absurdum*, I have given a simple proof that it is impossible to prove the RH by the mathematic way without first discovering the Prime Law.

So long the RH remains unproven, no one can claim that a near zero probability event could actually happen. But many have speculated it could and most scientists believe or speculate an arbitrary/random cause for the universe despite the near zero probability of such a randomly caused universe. To prove the RH by the mathematical way would only bolster such anti-reason beliefs. If there is a God, He/She would forbid that for sure. If you believe in God, you know you will be wasting time to prove the RH in the God-forbidden way.

Therefore, to prove the RH must entail eliminating any possibility of a lawless/arbitrary/random process in the creation of primes. The only way to do this is to discover the actual deterministic law of creating primes. So long you don't have the law, someone could always speculate a lawless process. If primes are arbitrarily created by a coin toss, the RH is much more likely to be false and may be easily disproved by finding a zero off the critical line. The fact that it has yet to be disproved suggests that a random process is unlikely. If a disproof of the RH indicates an arbitrary or lawless process in generating primes, a proof must exclude it. A proof simply cannot allow the possibility of a disproof. By allowing the possibility of an arbitrary/lawless/random process in creating primes, the possibility of a disproof shall always remain. If the way of seeking a proof is in principle incapable of excluding lawlessness in generating primes, it must have no chance of a successful proof.

Some mathematicians have doubted the importance of proving the RH in strictly mathematical terms, as Morris Kline (1908-1992) said in an interview: "If I could come back after five hundred years and find that the Riemann Hypothesis or Fermat's last 'theorem' was proved, I would be disappointed, because I would be pretty sure, in view of the history of attempts to prove these conjectures, that an enormous amount of time had been spent on proving theorems the are unimportant to the life of man."[21]. Indeed, a proof for the exact statement of the RH without relying on the Prime Law would most likely have little relevance to human life and would probably disappoint every mind. It would in fact have a negative impact on human value as it would support the absurd notion that an effectively zero probability event can actually happen. In contrast, proving the RH by the Prime Law provides a mathematical understanding of human creativity and affirms the human intuition that only laws can cause meaningful regularity.

There is another common sense proof why the RH cannot be proven via conventional mathematics. Doing so is to trivialize the RH. If the RH is the most fundamental unsolved mathematical question, which most think so, then it must be about the most fundamental/primary property of primes or the deepest mystery of primes.  The primary property in my view is the duality of uniqueness and uniformity and the Prime Law.  Common sense demands us to use primary property to deduce secondary properties which in turn deduce tertiary properties and so on.  We simply cannot deduce the primary property from secondary properties or deduce one secondary property from another parallel secondary property.  (The hardness of diamond is a secondary deduction of the carbon building blocks and their arrangements and cannot be deduced from another secondary property like the sparking property.)

      Mathematics including the RH are above the primary property of primes and are secondary or tertiary properties.  If the RH is fundamental, it must be either the primary property or an immediate deduction of the primary property.  It simply cannot be a tertiary property deduced from some secondary properties.  If it is secondary, then it also cannot be deduced from other parallel secondary properties.  If the RH can be proven using the secondary properties, it would be a tertiary property and as such it simply cannot be fundamental.  Here thus lies the absurdity.  Before proving the RH, it looks to be fundamental.  But after proving it, it becomes a trivial tertiary property.  Here, the primary questions about primes remain unsolved: is there a deterministic law of creating primes?  Is prime predictable?  What is the essence of primes?  Can that essence deduce all other secondary/tertiary properties including the RH?  The author Erica Klarreich put it: "Proving the Riemann Hypothesis won't end the story. It will prompt a sequence of even harder, more penetrating questions. Why do the primes achieve such a delicate balance between randomness and order?"[31].  It makes no sense for a proof of the RH, the most fundamental mathematical problem, to tell us nothing about the primary property of primes.  Only by proving the RH the correct way via the Prime Law, can we hope to end the story of primes.

      It is common sense that certain secondary property must be a direct deduction of the primary property and simply cannot be deduced any other way.  The RH is most likely such a secondary property.  Only by being so, can the RH claim to be fundamental next only to the primary essence of primes.  If not, we would have to find an intermediate secondary connection between the RH as the tertiary and the primary essence of primes.  That intermediate secondary connection would be more fundamental than the RH but is nowhere in sight in all of mathematics.  The only way to prove the RH while still maintain its status of being fundamental is to show that it is a direct/immediate deduction from the primary essence of primes.  The present notion of the essence of primes, i.e., indivisibility, is in fact a secondary property.  From that, one will never be able to reach the RH if the RH is fundamental.

Therefore, the great logician Kurt Godel is most likely correct in his view, as the mathematician Marcus du Sautoy wrote: "Godel himself had voiced such concerns in relation to the Riemann Hypothesis: perhaps the axioms that formed the foundations of the mathematical edifice were not broad enough to carry the required proof, in which case you might continue building upwards and never find a connection to the Hypothesis.  However, he did offer some consolation.  Godel believed that any conjecture of genuine interest cannot be forever out of reach.  It was just a matter of finding a new foundation stone to extend the base of the edifice.  Only by going back to the subject's foundations and seeking to broaden them would you be able to build up to the missing proof." [16].  The trinity of uniqueness, uniformity and the Prime Law is the new and broadened foundation for the primes.  As far as I know, the work here is the first attempt to prove the RH the Godel way "by going back to the subject's foundations."

This brings up the third way of proving that it is not possible for conventional mathematics to prove the RH.  Such a proof would trivialize any newly discovered foundations for the primes.  It would imply that the existing foundations of the primes are already complete.  But this notion can be easily proven to be false, because the existing definition of primes is merely a tautology.  Primes are not presently defined by their foundations (i.e., concepts more fundamental than numbers) but are circularly defined by concepts even less fundamental than numbers.  The ultimate understanding of the essence of primes should have a significant impact on mathematics.  It should solve some deepest puzzles of primes that should remain unsolvable otherwise.  For its impact to be negligible, which would be the case if the RH can be proven independently by conventional mathematics, is against reason and common sense, and in turn against how mathematics works.

I therefore conclude that the only way to prove the RH is by way of solving the ultimate question about primes, i.e., is there a deterministic law of creating primes that nonetheless allows primes to have the inherent property of unpredictability?  If the RH is directly related to the most fundamental essence of primes as many believe, then any method of proving the RH simply cannot succeed if it cannot prove the most fundamental concept about primes, i.e., an inherently unpredictable creature or prime is still a product of a fully deterministic or lawful process.

*The ultimate starting point of mathematics must be proven truth*

In mathematics, a proof is a demonstration that, assuming certain axioms, some statement is necessarily true.  A proof is derived by principles of deduction from certain axioms or common truth.  To prove the RH, one needs to be able to logically deduce it from some basic fundamental truth or axioms or essence of primes. In mathematics, axioms are neither derived by principles of deduction, nor are they demonstrable by formal proofs. Instead, an axiom is taken for granted as valid, and serves as a necessary starting point for deducing logically consistent propositions. But today's axioms can easily become tomorrow's deductions. From the essence of primes as uniqueness, we can deduce or prove that primes are indivisible. That prime is uniqueness is not an axiom but is a proven statement because the Prime Law can create primes by simply creating uniqueness.

Kurt Godel has proven that axioms drained of meaning and truth could never be used to prove that they will never lead to deductions that could contradict the axioms or other deductions. He further proved that any consistent axiom system free of human intuition is necessarily incomplete in that there will be true statements that can't be deduced from the axioms. Any mathematical edifice or formal system built on axioms drained of meaning and truth cannot have the certainty of truth and reality. The reality of mathematics is directly related to the meaning and truth of the axioms and in turn the intuition of the mind. But since not all intuitions are trustworthy, formalists like David Hilbert have tried to get rid of intuitions out of mathematics, which has been proven to be futile by Godel.

Here I suggest a novel approach to get rid of any uncertainty about the intuitions or intuition-based axioms. To have certainty and truth in the mathematical edifice, intuitions or axioms about the ultimate starting point need to be replaced by *proven* truth. Mathematics needs to be based on proven truth rather than axioms or intuitions. The ultimate law of the universe needs to have a firm foundation in proven truth rather than axioms. For the ultimate starting points to be proven truth rather than unproven axioms, they have to prove themselves or be a deduction of themselves. They must be self-caused.

The essence of prime is the starting point of primes, which has three parts, the duality of uniqueness and uniformity and the mind programmed with the Prime Law. The duality starting point is a proven truth because it has caused by way of the Prime Law an infinite number of itself, the primes, which are nothing but the duality. The Prime Law is also self-caused because it is caused or created by the human mind that is already programmed with the Prime Law. The creation law that programs the human mind must be the cause of every law or theory that is formulated by the mind, including the creation law or the Prime Law itself. In this sense, the Prime Law has caused itself. Therefore, the starting point of primes, the duality building block and the Prime Law, has caused or proved itself. Thus, mathematics built upon the essence of primes must be complete and consistent. Since the mind is part of the Prime Law, mathematics simply cannot be complete and consistent without the mind or the *proven* intuition of the mind. This independently confirms what Godel has found.

It is much more reassuring to know that the starting point of all reality is a proven position rather than an assumption or axiom.  If we only have an axiom, we will never know if it may be merely a deduction of a deeper truth.  Thus, if we can prove our starting point, which essentially means that the starting point can cause or deduce itself, it must mean that we have truly discovered the starting point.  Everything except the starting point must have a cause other than itself.

In the ideal world of logical and coherent knowledge, nothing is taken for granted without a proof, especially the starting point of everything. If a purported starting point is inherently un-provable, such as randomness, it is an axiom that may simply be an illusion.  Randomness as starting point also cannot qualify as a scientific position since it cannot be tested to be true or false.  It deduces all kinds of outcomes and excludes none.  A proof in the mathematical sense by using the principle of deduction is impossible with a starting point of randomness.  If a statement or theory cannot be proven by logical coherence in the mathematical sense of proof, it is not much of a theory by the standard of mathematics.  Given the existence of so many proven mathematical theorems and their 'unreasonable effectiveness' in understanding nature, it is impossible to imagine that the most fundamental law of nature should turn out to be not up to the standard of mathematics.

**Summary**

Many have observed that beauty lies in the seemingly impossible unity of extreme bipolar opposites.  The ultimate complexity should be linked with the ultimate simplicity.  The beauty of primes is the seemingly impossible bipolar unity of determinism and unpredictability. The only way to prove this beauty is to find the deterministic law that nonetheless produces outcomes with the inherent property of unpredictability.  To prove the RH by conventional mathematics simply cannot be done because it is equivalent to proving the impossible possible.

For the RH, a complex problem, to have a simple intuitive proof is the most beautiful thing one could hope for and is in full accord with how nature works, as nature is elegant, simple, and effortless. It would certainly be ugly and awkward for the RH to be proven by complex mathematics that only a handful of people can really understand or care to read.  The law behind human creativity is very simple, whose iterative use creates all complex creatures of the mind.  The lawfulness of the creative process makes the process sustainable and endless. It makes creation an achievable goal for any individual mind but an endless goal for the whole humanity. The unpredictability of the creative outcomes makes creation meaningful and challenging to the mind. The scientific understanding of nature is ultimately based on mathematics and creativity. The two are literally one and the same as they share a common foundation in the Prime Law. That human creativity is coded by numbers explains why nature can be understood by human creativity through mathematics.

**Acknowledgements:**

I am deeply grateful to Professor Louis H. Kauffman, Department of Mathematics, Statistics, and Computer Science at the University of Illinois at Chicago, for a detailed critic of the manuscript and for providing many helpful and insightful comments. I also thank 3 anonymous reviewers who have provided helpful comments on earlier versions of this manuscript. This work was supported by a grant from the NIH (RO1 CA 105347).

**Table 1. Creating primes.**

**A.** The contents of the reality domain at the time of creation. **B.** Subsequent progression of the reality domain. From left to right represents the number species with each number increasing in value from the previous number by the unit value of the beginning number; the species terminates at the Nth position where N is the numeric value of the last known prime (N>2). The species of 2 is listed not because 2 is a prime but because it is an inseparable part of creating the first prime 1. Successive prime numbers from small to large are listed on the left side column in the order from top to bottom. The Table can be expanded in a prime by prime manner over time to infinity, in both the vertical direction from top to bottom and the lateral direction from left to right.

A.

| | |
|---|---|
| 1 | 2 |
| 2 | 4 |

B.

| | | | | | | |
|---|---|---|---|---|---|---|
| 1 | 2 | 3 | 4 | 5 | ... | N |
| 2 | 4 | 6 | 8 | 10 | ... | 2N |
| 3 | 6 | 9 | 12 | 15 | ... | 3N |
| 5 | 10 | 15 | 20 | 25 | ... | 5N |
| ... | ... | ... | ... | ... | ... | ... |
| N | 2N | 3N | 4N | 5N | ... | NN |